\documentclass{elsart}
\usepackage{amssymb}
\usepackage{pstricks, pst-plot, pst-node, pst-tree}
\usepackage{array}
\usepackage{fancyhdr}
\usepackage{longtable}
\usepackage{tabularx}
\usepackage{multirow}
\usepackage{rotating}

\usepackage{makeidx}
\usepackage{amsmath}
\usepackage{theorem}
\usepackage[mathscr]{eucal}
\usepackage{enumerate}
\usepackage[dvips]{epsfig}
\usepackage{float}
\journal{Mathematics and Computers in Simulation}
\volume{77}
\pubyear{2008}
\newcounter{mylastpage}
\issue{4}
\usepackage{hyperref}
\makeatletter
\def\ps@copyright{%
 \let\@oddhead\@empty
 \let\@evenhead\@empty
 \def\@oddfoot{\small\slshape\hskip-7em
   Published in \@journal\ \@volume\ (\the\@pubyear)\ no.\ \@issue, pp.\ \ESpagenumber{firstpage}--\ESpagenumber{mylastpage},
\href{http://dx.doi.org/10.1016/j.matcom.2007.04.016}{doi: 10.1016/j.matcom.2007.04.016}}%
 \let\@evenfoot\@oddfoot}
\makeatother
\theoremstyle{plain}
\newtheorem{Def}{Definition}[section]

\newtheorem{The}[Def]{Theorem}

\newcommand{\arsinh}{\operatorname{arsinh}}

\newcommand{\E}{\operatorname{E}}
\newcommand{\Prob}{\operatorname{P}}
\begin{document}

\begin{frontmatter}

\title{Classification of Stochastic Runge--Kutta Methods
    for the Weak Approximation of Stochastic Differential Equations}
\author{Kristian Debrabant} and
\ead{debrabant@mathematik.tu-darmstadt.de}
\author{Andreas R\"{o}{\ss}ler}
\ead{roessler@mathematik.tu-darmstadt.de}
\address{Technische Universit\"{a}t Darmstadt, Fachbereich Mathematik, Schlo{\ss}gartenstr.7,
D-64289 Darmstadt, Germany}
\begin{abstract}
In the present paper, a class of stochastic Runge--Kutta methods
containing the second order stochastic Runge--Kutta scheme due
to E.~Platen for the weak approximation of It\^{o} stochastic
differential equation systems with a multi--dimensional Wiener
process is considered. Order one and order two conditions for the
coefficients of explicit stochastic Runge--Kutta methods are solved
and the solution space of the possible coefficients is analyzed. A
full classification of the coefficients for such stochastic
Runge--Kutta schemes of order one and two with minimal stage numbers
is calculated. Further, within the considered class of stochastic
Runge--Kutta schemes coefficients for optimal schemes in the sense
that additionally some higher order conditions are fulfilled are
presented.
\end{abstract}

\end{frontmatter}
\section{Introduction} \label{Introduction}
Recently, the development of numerical schemes for strong as well as
weak approximation of stochastic differential equations (SDEs) has
focused amongst others on Runge--Kutta type
schemes~\cite{BuBu00a,BuBu96,KP99,KoMiSu97,MacNav01,Mil95,Mil04,Roe06a,Roe06b,Roe03,ToVA02}.
This is due to the increasing complexity of stochastic Taylor
expansions and the desire to avoid derivatives in higher order
approximation schemes. In section~\ref{Sec:SRK-methods}, a class of
stochastic Runge--Kutta (SRK) methods due to
R\"{o}{\ss}ler~\cite{Roe06a,Roe06b,Roe03} for the weak approximation of It\^{o}
SDE systems with a multi--dimensional Wiener process is considered.
This class contains as a special case the second order SRK scheme
proposed by Platen~\cite{KP99} as well as the class of SRK methods
proposed by Tocino and Vigo-Aguiar~\cite{ToVA02}. Order conditions
for coefficients of the SRK methods have been calculated by applying
the multi-colored rooted tree analysis due to
R\"{o}{\ss}ler~\cite{Roe06a,Roe06b,Roe03}. In contrast to earlier work on
this topic, the aim of the present paper is to analyze these order
conditions with the objective to determine a full classification of
the coefficients for this class of SRK methods. A full
classification for order one SRK schemes with $s=1$ stage and for
order one SRK schemes with deterministic order two for $s=2$ stages
as well as for second order SRK schemes with $s=3$ stages is
calculated in section~\ref{Sec:Parameter-families}. Further, some
optimal schemes are derived from this classification in
section~\ref{Sec:Optimal-schemes} by taking into account additional
higher order conditions. Their performance is studied by some
numerical examples in section~\ref{Sec:Numerical-Example}.
\\ \\
We denote by $(X_t)_{t \in I}$ the solution of the $d$-dimensional
It\^{o} SDE defined by
\begin{equation} \label{St-lg-sde-ito-1}
    {\mathrm{d}} X_t = a(t,X_t) \, {\mathrm{d}}t + b(t,X_t) \,
    {\mathrm{d}}W_t, \qquad X_{t_0} = x_{0},
\end{equation}
with an $m$-dimensional Wiener process $(W_t)_{t \geq 0}$ and
$I=[t_0,T]$.
We assume that the Borel-measurable coefficients $a : I \times
\mathbb{R}^d \rightarrow \mathbb{R}^d$ and $b : I \times
\mathbb{R}^d \rightarrow \mathbb{R}^{d \times m}$ satisfy a
Lipschitz and a linear growth condition such that the Existence
and Uniqueness Theorem~\cite{KP99} applies. In the following, let
$b^j(t,x) = (b^{i,j}(t,x))_{1 \leq i \leq d} \in \mathbb{R}^d$
denote the $j$th column of the diffusion matrix $b(t,x)$ for $j=1,
\ldots, m$. \\ \\
Let a discretization $I_h = \{t_0, t_1, \ldots, t_N\}$ with $t_0 <
t_1 < \ldots < t_N =T$ of the time interval $I=[t_0,T]$ with step
sizes $h_n = t_{n+1}-t_n$ for $n=0,1, \ldots, N-1$ be given.
Further, define $h = \max_{0 \leq n < N} h_n$ as the maximum step
size.
Let $C_P^l(\mathbb{R}^d, \mathbb{R})$ denote the space of all $g
\in C^l(\mathbb{R}^d,\mathbb{R})$ fulfilling a polynomial growth
condition and let $g \in C_P^{k,l}(I \times \mathbb{R}^d,
\mathbb{R})$ if $g(\cdot,x) \in C^{k}(I,\mathbb{R})$ and
$g(t,\cdot) \in C_P^l(\mathbb{R}^d, \mathbb{R})$
for all $t \in I$ and $x \in \mathbb{R}^d$ \cite{KP99}.
\begin{Def}
    A time discrete approximation
    $Y=(Y_t)_{t \in I_h}$
    converges weakly
    with order $p$ to $X$ as $h \rightarrow 0$ at time $t \in I_h$ if
    for each $f \in C_P^{2(p+1)}(\mathbb{R}^d, \mathbb{R})$
    exist a constant $C_f$
    and a finite $\delta_0 > 0$ such that
    \begin{equation}
        | \E(f(X_t)) - \E(f(Y_t)) | \leq C_f \, h^p
    \end{equation}
    holds for each $h \in \, ]0,\delta_0[\,$.
\end{Def}
\section{Stochastic Runge--Kutta Methods}
\label{Sec:SRK-methods}
We consider stochastic Runge--Kutta methods as proposed in
\cite{Roe06a,Roe06b,Roe03} for the weak approximation of
SDE~(\ref{St-lg-sde-ito-1}). Therefore, the $d$-dimensional
approximation process $Y$ of an explicit $s$-stage SRK method is
defined by $Y_{t_0} = x_0$ and
\begin{equation} \label{SRK-method-Ito-Wm-allg01}
    \begin{split}
    Y_{t_{n+1}} = Y_{t_n} & + \sum_{i=1}^s
    \alpha_i \, a(t_n+c_i^{(0)} h_n, H_i^{(0)}) \, h_n \\
    & + \sum_{i=1}^s
    \sum_{k=1}^m
    \beta_i^{(1)} \, b^{k}(t_n+c_i^{(1)} h_n, H_i^{(k)}) \, \hat{I}_{(k)} \\
    & +
    \sum_{i=1}^s \sum_{k=1}^m
    \beta_i^{(2)} \, b^{k}(t_n+c_i^{(1)} h_n, H_i^{(k)}) \,
    \tfrac{\hat{I}_{(k,k)}}{\sqrt{h_n}} \\
    & + \sum_{i=1}^s
    \sum_{\substack{k,l=1 \\ k \neq l}}^m
    \beta_i^{(3)} \, b^{k}(t_n+c_i^{(2)} h_n, \hat{H}_i^{(l)}) \, \hat{I}_{(k)} \\
    & +
    \sum_{i=1}^s \sum_{\substack{k,l=1 \\ k \neq l}}^m
    \beta_i^{(4)} \, b^{k}(t_n+c_i^{(2)} h_n, \hat{H}_i^{(l)}) \,
    \tfrac{\hat{I}_{(k,l)}}{\sqrt{h_n}}
    \end{split}
\end{equation}
for $n=0,1, \ldots, N-1$ with stage values
\begin{alignat*}{5}
    H_i^{(0)} &=&\,\, Y_{t_n} &+ \sum_{j=1}^{i-1} A_{ij}^{(0)}
    \, a(t_n+c_j^{(0)} h_n, H_j^{(0)}) \, h_n \\
    && &+ \sum_{j=1}^{i-1} \sum_{r=1}^m
    B_{ij}^{(0)} \, b^r(t_n+c_j^{(1)} h_n, H_j^{(r)}) \, \hat{I}_{(r)} \\
    H_i^{(k)} &=&\,\, Y_{t_n} &+ \sum_{j=1}^{i-1} A_{ij}^{(1)}
    \, a(t_n+c_j^{(0)} h_n, H_j^{(0)}) \, h_n \\
    && &+ \sum_{j=1}^{i-1}
    B_{ij}^{(1)} \, b^k(t_n+c_j^{(1)} h_n, H_j^{(k)}) \,
    \sqrt{h_n} \\
    \hat{H}_i^{(k)} &=&\,\, Y_{t_n} &+ \sum_{j=1}^{i-1} A_{ij}^{(2)}
    \, a(t_n+c_j^{(0)} h_n, H_j^{(0)}) \, h_n \\
    && &+ \sum_{j=1}^{i-1}
    B_{ij}^{(2)} \, b^k(t_n+c_j^{(1)} h_n, H_j^{(k)}) \, \sqrt{h_n}
\end{alignat*}
for $i=1, \ldots, s$ and $k=1, \ldots, m$. Here, $\alpha,
\beta^{(1)},\dots,\beta^{(4)},c^{(q)}\in \mathbb{R}^s$ and
$A^{(q)}$, $B^{(q)} \in \mathbb{R}^{s \times s}$ for $0\leq q \leq
2$ with $A_{ij}^{(q)} = B_{ij}^{(q)} = 0$ for $j \geq i$ are the
vectors and matrices of coefficients of the SRK method. We choose
$c^{(q)}=A^{(q)} e$ for $0 \leq q \leq 2$ with a vector $e=(1,
\ldots, 1)^T$ \cite{Roe06a}. In the following, the product of column
vectors is defined component-wise. The coefficients of the SRK
method~(\ref{SRK-method-Ito-Wm-allg01}) are determined by the
following Butcher tableau:
\renewcommand{\arraystretch}{1.8}
\begin{center}
\begin{tabular}{c|c|c|c}
    $c^{(0)}$ & ${A}^{(0)}$ & $B^{(0)}$ & \\
    \cline{1-4}
    $c^{(1)}$ & ${A}^{(1)}$ & $B^{(1)}$ & \\
    \cline{1-4}
    $c^{(2)}$ & ${A}^{(2)}$ & $B^{(2)}$ & \\
    \hline
    & $\alpha^T$ & ${\beta^{(1)}}^T$ & ${\beta^{(2)}}^T$ \\
    \cline{2-4}
    & & ${\beta^{(3)}}^T$ & ${\beta^{(4)}}^T$
\end{tabular}
\end{center}
\renewcommand{\arraystretch}{1.0}
The random variables of the SRK method are defined by three-point
distributed random variables with $\Prob(\hat{I}_{(r)} = \pm
\sqrt{3 \, h_n} ) = \frac{1}{6}$ and $\Prob(\hat{I}_{(r)} = 0 ) =
\frac{2}{3}$. Further, $ \hat{I}_{(k,l)} = \frac{1}{2} (
\hat{I}_{(k)} \, \hat{I}_{(l)} + V^{k,l} )$. The $V^{k,l}$ are
independent two-point distributed random variables with
$\Prob(V^{k,l} = \pm h_n) = \frac{1}{2}$ for $l=1, \ldots, k-1$,
$V^{k,k} = -h_n$ and $V^{k,l} = -V^{l,k}$ for $l=k + 1, \ldots, m$
and $k=1, \ldots, m$ \cite{KP99}.
\\ \\
By the application of the multi--colored rooted tree
analysis~\cite{Roe06a,Roe03}, order conditions for the
coefficients of the SRK method~(\ref{SRK-method-Ito-Wm-allg01})
can be easily determined.
As a result of this, the following
Theorem~\ref{SRK-theorem-ito-ord2-Wm-main1} due to
R\"{o}{\ss}ler~\cite{Roe03} gives order conditions for the SRK
method~(\ref{SRK-method-Ito-Wm-allg01}) up to order two.
\begin{The} \label{SRK-theorem-ito-ord2-Wm-main1}
    Let $a^i, b^{ij} \in C_P^{2,4}(I \times
    \mathbb{R}^d,\mathbb{R})$ for
    $1 \leq i \leq d$, $1 \leq j \leq m$.
    If the coefficients of the SRK
    method~(\ref{SRK-method-Ito-Wm-allg01})
    fulfill the equations
    \begin{alignat*}{5}
        1&. \quad \alpha^T e = 1 \qquad \qquad
        &2. \quad &{\beta^{(4)}}^T e = 0 \qquad \qquad
        &3. \quad &{\beta^{(3)}}^T e = 0 \\
        4&. \quad ({\beta^{(1)}}^T e)^2 = 1 \qquad \qquad
        &5. \quad &{\beta^{(2)}}^T e = 0 \qquad \qquad
        &6. \quad &{\beta^{(1)}}^T {B^{(1)}} e = 0 \\
        7&. \quad {\beta^{(3)}}^T
        {B^{(2)}} e = 0
    \end{alignat*}
    then the SRK method converges with order 1 in the weak sense.
    In addition,
    if $a^i, b^{ij} \in C_P^{3,6}(I \times
    \mathbb{R}^d,\mathbb{R})$ for
    $1 \leq i \leq d$, $1 \leq j \leq m$
    and if the equations
    {\allowdisplaybreaks
    \begin{alignat*}{3}
        8&. \quad \alpha^T A^{(0)} e = \tfrac{1}{2}
        \qquad \qquad
        &9. \quad &\alpha^T (B^{(0)} e)^2 = \tfrac{1}{2} \\
        10&. \quad ({\beta^{(1)}}^T e) (\alpha^T B^{(0)} e) =
        \tfrac{1}{2}
        \qquad \qquad
        &11. \quad &({\beta^{(1)}}^T e) ({\beta^{(1)}}^T A^{(1)} e) = \tfrac{1}{2} \\
        12&. \quad {\beta^{(3)}}^T A^{(2)} e = 0
        \qquad \qquad
        &13. \quad &{\beta^{(2)}}^T B^{(1)} e = 1 \\
        14&. \quad {\beta^{(4)}}^T B^{(2)} e = 1
        \qquad \qquad
        &15. \quad &({\beta^{(1)}}^T e) ({\beta^{(1)}}^T
        (B^{(1)} e)^2) = \tfrac{1}{2} \\
        16&. \quad ({\beta^{(1)}}^T e) ({\beta^{(3)}}^T
        (B^{(2)} e)^2) = \tfrac{1}{2} \qquad
        &17. \quad &{\beta^{(1)}}^T (B^{(1)} (B^{(1)} e)) =
        0 \\
        18&. \quad {\beta^{(3)}}^T (B^{(2)}
        (B^{(1)}e)) = 0 \qquad
        &19. \quad &{\beta^{(3)}}^T (A^{(2)} (B^{(0)}
        e)) = 0 \\
        20&. \quad {\beta^{(1)}}^T (A^{(1)} (B^{(0)} e)) =
        0 \qquad
        &21. \quad &\alpha^T (B^{(0)} (B^{(1)} e)) = 0 \\
        22&. \quad {\beta^{(2)}}^T A^{(1)} e = 0 \qquad
        &23. \quad &{\beta^{(4)}}^T A^{(2)} e = 0 \\
        24&. \quad {\beta^{(1)}}^T ((A^{(1)} e) (B^{(1)} e)) = 0
        \qquad
        &25. \quad &{\beta^{(3)}}^T ((A^{(2)} e)
        (B^{(2)} e)) = 0 \\
        26&. \quad {\beta^{(4)}}^T (A^{(2)} (B^{(0)}
        e)) = 0 \qquad
        &27. \quad &{\beta^{(2)}}^T (A^{(1)} (B^{(0)} e))
        = 0 \\
        28&. \quad {\beta^{(2)}}^T (A^{(1)} (B^{(0)} e)^2) = 0
        \qquad
        &29. \quad &{\beta^{(4)}}^T (A^{(2)} (B^{(0)}
        e)^2) = 0 \\
        30&. \quad {\beta^{(3)}}^T (B^{(2)} (A^{(1)}
        e)) = 0 \qquad
        &31. \quad &{\beta^{(1)}}^T (B^{(1)} ( A^{(1)} e)) =
        0 \\
        32&. \quad {\beta^{(2)}}^T (B^{(1)} e)^2 = 0 \qquad
        &33. \quad &{\beta^{(4)}}^T (B^{(2)} e)^2 = 0 \\
        34&. \quad {\beta^{(4)}}^T (B^{(2)} (B^{(1)}
        e)) = 0 \qquad
        &35. \quad &{\beta^{(2)}}^T (B^{(1)} (B^{(1)} e))
        = 0 \\
        36&. \quad {\beta^{(1)}}^T (B^{(1)} e)^3 = 0 \qquad
        &37. \quad &{\beta^{(3)}}^T (B^{(2)} e)^3 = 0 \\
        38&. \quad {\beta^{(1)}}^T (B^{(1)} (B^{(1)}
        e)^2) = 0 \qquad
        &39. \quad &{\beta^{(3)}}^T (B^{(2)} (B^{(1)}
        e)^2) = 0
    \end{alignat*}
}
    \begin{alignat*}{3}
        40&. \quad \alpha^T ((B^{(0)} e) (B^{(0)}
        (B^{(1)} e))) = 0 \\
        41&. \quad {\beta^{(1)}}^T ((A^{(1)} (B^{(0)}
        e)) (B^{(1)} e)) = 0 \\
        42&. \quad {\beta^{(3)}}^T ((A^{(2)}
        (B^{(0)} e)) (B^{(2)} e)) = 0 \\
        43&. \quad {\beta^{(1)}}^T (A^{(1)} (B^{(0)} (B^{(1)} e))) = 0 \\
        44&. \quad {\beta^{(3)}}^T (A^{(2)} (B^{(0)}
        (B^{(1)} e))) = 0 \\
        45&. \quad {\beta^{(1)}}^T (B^{(1)} (A^{(1)}
        (B^{(0)} e))) = 0 \\
        46&. \quad {\beta^{(3)}}^T (B^{(2)} (A^{(1)}
        (B^{(0)} e))) = 0 \\
        47&. \quad {\beta^{(1)}}^T ((B^{(1)} e) (B^{(1)}
        (B^{(1)} e))) = 0 \\
        48&. \quad {\beta^{(3)}}^T ((B^{(2)} e)
        (B^{(2)} (B^{(1)} e))) = 0 \\
        49&. \quad {\beta^{(1)}}^T (B^{(1)} (B^{(1)}
        (B^{(1)} e))) = 0 \\
        50&. \quad {\beta^{(3)}}^T (B^{(2)} (B^{(1)}
        (B^{(1)} e))) = 0
    \end{alignat*}
    are fulfilled then the stochastic Runge--Kutta
    method~(\ref{SRK-method-Ito-Wm-allg01}) converges with order~2
    in the weak sense.
\end{The}
In the case of $m>1$ one has to solve 50 non-linear equations in
order to calculate coefficients for an order two SRK
method~(\ref{SRK-method-Ito-Wm-allg01}). However, in the case of
$m=1$ these conditions are reduced to 28 equations which have to be
solved~\cite{Roe06b}. Thus, the analysis of the space of all
admissible coefficients is not an easy task. It turns out that
explicit order one SRK methods need at least $s=1$ stage while order
two SRK methods need $s \geq 3$ stages. This is due to  e.g.\
conditions 6.\ and 15., which can not be fulfilled in the case of
$s \leq 2$ stages for explicit order two SRK methods.
In the following, we distinguish between the stochastic and the
deterministic order of convergence. Let $p_S=p$ denote the order of
convergence of the SRK method if it is applied to an SDE and let
$p_D$ with $p_D \geq p_S$ denote the order of convergence of the SRK
method if it is applied to a deterministic ordinary differential
equation (ODE), i.e., SDE~(\ref{St-lg-sde-ito-1}) with $b \equiv 0$.
We also write $(p_D,p_S)$ in the following \cite{Roe06b,Roe03}.
\section{Parameter Families for SRK Methods}
\label{Sec:Parameter-families}
\subsection{Coefficients for SRK Methods of Order (1,1)}
\label{Sec:Coeff-SRK-Ord-1-1}
First, we analyze explicit SRK
methods~(\ref{SRK-method-Ito-Wm-allg01}) of order $p_D=p_S=1$ with
$s=1$ stage. Considering the order one conditions 1.--7.\ in
Theorem~\ref{SRK-theorem-ito-ord2-Wm-main1}, the corresponding
coefficients are uniquely determined for $c_1 \in \{-1,1\}$ by
\begin{equation} \label{Parameter-Ord11-all}
    \alpha_1 = 1, \qquad \beta_1^{(1)} = c_1,
    \qquad \beta_1^{(2)} = 0, \qquad
    \beta_1^{(3)} = 0, \qquad
    \beta_1^{(4)} = 0.
\end{equation}
The resulting class of SRK schemes coincides with the well-known
Euler-Maruyama scheme.
\subsection{Coefficients for SRK Methods of Order (2,1)}
\label{Sec:Coeff-SRK-Ord-2-1}
Next, we consider the case of $s=2$ stage explicit SRK methods
(\ref{SRK-method-Ito-Wm-allg01}). As already mentioned in
Section~\ref{Sec:SRK-methods}, it is not possible to attain order
$p_S=2$. However, we can find some SRK methods of order $p_D=2$ and
$p_S=1$ corresponding to the following parameter family: From
condition 1.\ of Theorem~\ref{SRK-theorem-ito-ord2-Wm-main1} follows
$\alpha_1 = 1 -\alpha_2$ and taking into account the order 2
condition 8.\ we obtain $\alpha_2 = \frac{1}{2 A^{(0)}_{21}}$ for
$A^{(0)}_{21} \neq 0$. Further, condition 2.\ yields $\beta^{(4)}_1
= -\beta^{(4)}_2$, condition 3.\ results in $\beta^{(3)}_1 =
-\beta^{(3)}_2$ and condition 5.\ is fulfilled if $\beta^{(2)}_1 =
-\beta^{(2)}_2$ while condition 4.\ holds for $\beta^{(1)}_1 = c_1 -
\beta^{(1)}_2$ with $c_1 \in \{-1,1\}$. Finally, considering
condition 6.\ we need that $\beta^{(1)}_2 = 0$ or $B^{(1)}_{21}=0$
and for condition 7.\ analogously that $\beta^{(3)}_2 = 0$ or
$B^{(2)}_{21}=0$ hold. Thus, this class of SRK methods is determined
by
\begin{alignat}{7}
    \alpha^T &= \begin{bmatrix} 1-\frac{1}{2 c_2} && \frac{1}{2
    c_2} \end{bmatrix} , &\quad \quad
    {\beta^{(1)}}^T &= \begin{bmatrix} c_1-c_4 && c_4 \end{bmatrix},
    &\quad \quad
    {\beta^{(2)}}^T &= \begin{bmatrix} c_5 && -c_5 \end{bmatrix}, \notag \\
    && {\beta^{(3)}}^T &= \begin{bmatrix} c_6 && -c_6 \end{bmatrix},
    &\quad
    {\beta^{(4)}}^T &= \begin{bmatrix} c_7 && -c_7 \end{bmatrix}, \notag \\
    A^{(0)} &= \begin{bmatrix} 0 && 0 \\ c_2 && 0 \end{bmatrix}, &\quad
    A^{(1)} &= \begin{bmatrix} 0 && 0 \\ c_8 && 0 \end{bmatrix}, &\quad
    A^{(2)} &= \begin{bmatrix} 0 && 0 \\ c_{9} && 0 \end{bmatrix},
    \notag \\
    B^{(0)} &= \begin{bmatrix} 0 && 0 \\ c_3 && 0 \end{bmatrix},
    &\quad
    B^{(1)} &= \begin{bmatrix} 0 && 0 \\ c_{10} && 0 \end{bmatrix}, &\quad
    B^{(2)} &= \begin{bmatrix} 0 && 0 \\ c_{11} && 0 \label{Parameter-Ord21-all}
    \end{bmatrix},
\end{alignat}
for $c_1 \in \{-1,1\}$ and $c_2, \ldots, c_{11} \in \mathbb{R}$
with $c_2 \neq 0$, $c_4 \, c_{10} =0$ and $c_6 \, c_{11} = 0$.
\subsection{Coefficients for SRK Methods of Order (2,2)}
\label{Sec:Coeff-SRK-Ord-2-2}
Now, we consider explicit SRK
methods~(\ref{SRK-method-Ito-Wm-allg01}) of order $p_D=p_S=2$ with
$s=3$ stages. Then, the SRK schemes of the class under consideration
are completely characterized by the following families of
coefficients which follow from the order conditions in
Theorem~\ref{SRK-theorem-ito-ord2-Wm-main1}: Due to conditions 13.\
and 32.\ we need $B^{(1)}_{21} \neq 0$ and from conditions 15.\ and
36.\ follows $\beta^{(1)}_3 \neq 0$. Thus, there exist no SRK
schemes of the considered class attaining order $p_D=p_S=2$ with
less than 3 stages. Now, by condition 17.\ follows that
$B^{(1)}_{32}=0$ and we deduce from 6., 15.\ and 36.\ that
$B^{(1)}_{31} = -B^{(1)}_{21}\neq0$. Analyzing the weights, we calculate
from conditions 5., 13.\ and 32.\ that $\beta^{(2)}_2 = \tfrac{1}{2
B^{(1)}_{21}}$, $\beta^{(2)}_3 = -\tfrac{1}{2 B^{(1)}_{21}}$ and
$\beta^{(2)}_1 = 0$. For $c_1 \in \{-1,1\}$ we obtain from conditions
4., 6.\ and 15.\ the weights $\beta^{(1)}_1 = c_1 - \tfrac{c_1}{2
(B^{(1)}_{21})^2}$ and $\beta^{(1)}_2 = \beta^{(1)}_3 =
\tfrac{c_1}{4 (B^{(1)}_{21})^2}$. Now, due to 24.\ and 11.\ we need
that $A^{(1)}_{21} = (B^{(1)}_{21})^2$ and $A^{(1)}_{31} =
(B^{(1)}_{21})^2 - A^{(1)}_{32}$. Applying now conditions 3., 7.,
16.\ and 37.\ we conclude that $B^{(2)}_{21} + B^{(2)}_{31} +
B^{(2)}_{32} = 0$, $B^{(2)}_{21} \neq 0$, $B^{(2)}_{21} \neq
B^{(2)}_{31} + B^{(2)}_{32}$ and that $\beta^{(3)}_1 =
-\tfrac{c_1}{2 (B^{(2)}_{21})^2}$, $\beta^{(3)}_2 = \tfrac{c_1}{4
(B^{(2)}_{21})^2}$ and $\beta^{(3)}_3 = \tfrac{c_1}{4
(B^{(2)}_{21})^2}$. Further, we can now determine the remaining
weights as $\beta^{(4)}_1 = 0$, $\beta^{(4)}_2 = \tfrac{1}{2
B^{(2)}_{21}}$ and $\beta^{(4)}_3 = -\tfrac{1}{2 B^{(2)}_{21}}$ from
conditions 2., 14.\ and 33., and we have $\alpha_1 = 1 - \alpha_2 -
\alpha_3$ due to condition 1. Then, we can consider condition 18.\
which needs $B^{(2)}_{32} = 0$ and we thus get with the previous
considerations that finally $B^{(2)}_{21} = -B^{(2)}_{31}$ has to be
fulfilled. Now, we obtain from conditions 12.\ and 23.\ that
$A^{(2)}_{21} = 0$ and that $A^{(2)}_{32} = -A^{(2)}_{31}$ has to be
fulfilled. Continuing in this manner, we have to distinguish the
following cases:
\begin{enumerate}[(A)]
    \item \label{Parameter-Fall1} For $\alpha_3 = 0$, the parameter
    family is
    given by $\alpha_1 = \alpha_2 = \tfrac{1}{2}$ and $B^{(0)}_{21}
    = c_1$, which follows from conditions 1., 9.\ and 10. Further,
    we calculate from condition 8.\ that $A^{(0)}_{21} = 1$, from
    20.\ that $A^{(1)}_{32} = 0$ and from condition 19.\ that
    $A^{(2)}_{32} = 0$.
    \item \label{Parameter-Fall2} For $\alpha_3 \neq 0$,
    condition 21.\ yields now that $B^{(0)}_{32} = 0$ and we have to
    consider the following cases:
    \begin{enumerate}
        \item \label{Parameter-Fall2.1} For $B^{(0)}_{21} = 0$
        it follows from conditions 9.\ and 10.\ that $B^{(0)}_{31} = c_1$ and $\alpha_3 = \tfrac{1}{2}$. Thus,
        by condition 1.\ it follows immediately that $\alpha_1 =
        \tfrac{1}{2} - \alpha_2$.
        \begin{enumerate}
            \item \label{Parameter-Fall2.1.1} If $A^{(0)}_{21} =
            0$ then condition 8.\ implies that in this case
            $A^{(0)}_{32} = 1 - A^{(0)}_{31}$ has to be fulfilled.
            \item \label{Parameter-Fall2.1.2} If $A^{(0)}_{21}
            \neq 0$ then condition 8.\ yields that $\alpha_1 =
            \tfrac{1}{2} - \tfrac{1-A^{(0)}_{31}-A^{(0)}_{32}}{2
            A^{(0)}_{21}}$ and $\alpha_2 = \tfrac{1-A^{(0)}_{31} -
            A^{(0)}_{32}}{2 A^{(0)}_{21}}$.
        \end{enumerate}
        \item \label{Parameter-Fall2.2} For $B^{(0)}_{21}
        \neq 0$, we calculate from condition 20.\ that $A^{(1)}_{32}
        = 0$ and from condition 19.\ that $A^{(2)}_{32} = 0$ which
        implies that also $A^{(2)}_{31} = 0$ due to $A^{(2)}_{32} =
        -A^{(2)}_{31}$. Now, with $\kappa = \alpha_2 \, \alpha_3 (2
        \alpha_2 + 2 \alpha_3 - 1)$ and conditions 9.\ and 10.\ we
        have to consider the following cases:
        \begin{enumerate}
            \item \label{Parameter-Fall2.2.1} In the case of
            $\kappa \geq 0$, $\alpha_2 \neq -\alpha_3$, $\alpha_2
            \neq 0$ and $\alpha_2 \neq \pm \sqrt{\kappa}$ it follows
            that $B^{(0)}_{21} = \tfrac{c_1 (\alpha_2 \mp \sqrt{\kappa})}{2
            \alpha_2 (\alpha_2+\alpha_3)}$ and $B^{(0)}_{31} =
            \tfrac{c_1 (\alpha_3 \pm \sqrt{\kappa})}{2 \alpha_3
            (\alpha_2 + \alpha_3)}$.
            \item \label{Parameter-Fall2.2.2} If $\alpha_2 =
            0$ and $\alpha_3 = \tfrac{1}{2}$ then we can conclude
            that $B^{(0)}_{31} = c_1$ has to hold.
            \item \label{Parameter-Fall2.2.3} For $\alpha_2 =
            -\alpha_3$ and $\alpha_2 \neq -\tfrac{1}{2}$ it follows that
            $B^{(0)}_{21} = c_1 ( \tfrac{1}{2}+\tfrac{1}{4 \alpha_2})$ and $B^{(0)}_{31} = c_1 (\tfrac{1}{2} -
            \tfrac{1}{4 \alpha_2})$ has to be fulfilled.
        \end{enumerate}
        Due to condition 8.\ it follows that
        $A^{(0)}_{31} = \tfrac{1 - 2 \alpha_2 \, A^{(0)}_{21}}{2
        \alpha_3} - A^{(0)}_{32}$.
    \end{enumerate}
\end{enumerate}
Finally, one can easily check that all remaining conditions of
Theorem~\ref{SRK-theorem-ito-ord2-Wm-main1}, which have not been
mentioned explicitly in our analysis, are fulfilled by each
parameter family and thus do not contribute any further restrictions
for the coefficients. \\ \\
Summarizing our results, we have the following classification for
the SRK schemes of order $p_D=p_S=2$ for the considered class with
$s=3$ stages: For $c_1 \in \{-1,1\}$ and $c_2, c_3, c_4, c_5 \in
\mathbb{R}$ with $c_3 \neq 0$ and $c_4 \neq 0$ it holds
\begin{alignat}{5}
    {\beta^{(1)}}^T &= \begin{bmatrix} c_1-\frac{c_1}{2 c_3^2} & &
    \frac{c_1}{4 c_3^2} & & \frac{c_1}{4 c_3^2} \end{bmatrix},
    &\quad \quad
    {\beta^{(2)}}^T &= \begin{bmatrix} 0 & & \frac{1}{2 c_3} & &
    -\frac{1}{2 c_3} \end{bmatrix}, \\
    {\beta^{(3)}}^T &= \begin{bmatrix} -\frac{c_1}{2 c_4^2} & &
    \frac{c_1}{4 c_4^2} & & \frac{c_1}{4 c_4^2} \end{bmatrix},
    &\quad \quad
    {\beta^{(4)}}^T &= \begin{bmatrix} 0 & & \frac{1}{2 c_4} & &
    -\frac{1}{2 c_4} \end{bmatrix}, \\
    A^{(1)} &= \begin{bmatrix} 0 && 0 && 0 \\ c_3^2 && 0 && 0 \\
    c_3^2-c_2 && c_2 && 0 \end{bmatrix}, &\quad \quad B^{(1)} &=
    \begin{bmatrix} 0 && 0 && 0 \\ c_3 && 0 && 0 \\ -c_3 && 0 && 0
    \end{bmatrix}, \label{Parameter-Ord22-A1-B1} \\
    A^{(2)} &= \begin{bmatrix} 0 && 0 && 0 \\ 0 && 0 && 0 \\
    c_5 && -c_5 && 0 \end{bmatrix}, &\quad \quad B^{(2)} &=
    \begin{bmatrix} 0 && 0 && 0 \\ c_4 && 0 && 0 \\ -c_4 && 0 && 0
    \end{bmatrix} . \label{Parameter-Ord22-A2-B2}
\end{alignat}
Now, the following cases are possible: \\ \\
In the case (\ref{Parameter-Fall1}), we get with $c_2=c_5=0$ in
(\ref{Parameter-Ord22-A1-B1})--(\ref{Parameter-Ord22-A2-B2}) that
\begin{equation} \label{Parameter-Coeff-1}
    \alpha^T = \begin{bmatrix} \frac12 & &
    \frac12 & & 0 \end{bmatrix} ,
    \quad
    A^{(0)} = \begin{bmatrix} 0 && 0 && 0 \\ 1 && 0 && 0 \\
     0&& 0 && 0 \end{bmatrix}, \quad B^{(0)} =
    \begin{bmatrix} 0 && 0 && 0 \\ c_1 && 0 && 0 \\0  &&0  && 0
    \end{bmatrix} ,
\end{equation}
with $A^{(0)}_{31}=A^{(0)}_{32}=B^{(0)}_{31}=B^{(0)}_{32}=0$
because these coefficients are not relevant for the scheme due to
$\alpha_3=0$. \\ \\
For the case (\ref{Parameter-Fall2.1.1}) we get with $c_6, c_7 \in
\mathbb{R}$ the coefficients
\begin{equation} \label{Parameter-Coeff-2.1.1}
    \alpha^T = \begin{bmatrix} \frac{1}{2}-c_6 &&
    c_6 && \frac{1}{2} \end{bmatrix},
    \quad
    A^{(0)} = \begin{bmatrix} 0 && 0 && 0 \\ 0 && 0 && 0 \\
    c_7 && 1-c_7 && 0 \end{bmatrix}, \quad B^{(0)} =
    \begin{bmatrix} 0 && 0 && 0 \\ 0 && 0 && 0 \\ c_1 && 0 && 0
    \end{bmatrix} .
\end{equation}
\quad \\
Considering the case (\ref{Parameter-Fall2.1.2}) we obtain for
$c_6, c_7, c_8 \in \mathbb{R}$ with $c_6 \neq 0$ that
\begin{equation} \label{Parameter-Coeff-2.1.2}
    \alpha^T = \begin{bmatrix} \frac{1}{2}-\frac{1-c_7-c_8}{2 c_6}
    && \frac{1-c_7-c_8}{2 c_6} && \frac{1}{2} \end{bmatrix},
    \,\,
    A^{(0)} = \begin{bmatrix} 0 && 0 && 0 \\ c_6 && 0 && 0 \\
    c_7 && c_8 && 0 \end{bmatrix},
    \,\,
    B^{(0)} =
    \begin{bmatrix} 0 && 0 && 0 \\ 0 && 0 && 0 \\ c_1 && 0 && 0
    \end{bmatrix} .
\end{equation}
\quad \\
Next, we have the case (\ref{Parameter-Fall2.2.1}) with $c_2=c_5=0$
in (\ref{Parameter-Ord22-A1-B1})--(\ref{Parameter-Ord22-A2-B2}),
$c_6, c_7, c_8,c_9 \in \mathbb{R}$ and with $c_6 \neq 0$ and $c_6 \neq
-c_7 \neq 0$. Then, it holds with $\kappa = c_6 c_7 (2 c_6 + 2 c_7
-1)$ and $\lambda = \frac{1-2 c_6 c_8}{2 c_7}$ for $c_6 \neq \pm
\sqrt{\kappa}$ and $\kappa \geq 0$ that
\begin{equation} \label{Parameter-Coeff-2.2.1}
    \alpha^T = \begin{bmatrix} 1-c_6-c_7 && c_6 && c_7
    \end{bmatrix}, \,\,
    A^{(0)} = \begin{bmatrix} 0 && 0 && 0 \\ c_8 && 0 && 0 \\
    \lambda-c_9 && c_9 && 0 \end{bmatrix}, \,\, B^{(0)} =
    \begin{bmatrix} 0 && 0 && 0 \\ \frac{c_1}{2}
    \frac{c_6 \mp \sqrt{\kappa}}{c_6 (c_6 + c_7)} && 0 && 0 \\
    \frac{c_1}{2} \frac{c_7 \pm \sqrt{\kappa}}{c_7 (c_6 + c_7)} && 0 && 0
    \end{bmatrix} .
\end{equation}
\quad \\
The case (\ref{Parameter-Fall2.2.2}) yields for $c_6, c_7, c_8 \in
\mathbb{R}$ with $c_8 \neq 0$ and $c_2=c_5=0$ in
(\ref{Parameter-Ord22-A1-B1})--(\ref{Parameter-Ord22-A2-B2}) the
coefficients
\begin{equation} \label{Parameter-Coeff-2.2.2}
    \alpha^T = \begin{bmatrix} \frac{1}{2}
    && 0 && \frac{1}{2} \end{bmatrix},
    \quad
    A^{(0)} = \begin{bmatrix} 0 && 0 && 0 \\ c_6 && 0 && 0 \\
    1-c_7 && c_7 && 0 \end{bmatrix}, \quad B^{(0)} =
    \begin{bmatrix} 0 && 0 && 0 \\ c_8 && 0 && 0 \\ c_1 && 0 && 0
    \end{bmatrix} .
\end{equation}
\quad \\
Finally, we have the case (\ref{Parameter-Fall2.2.3}) for $c_6, c_7,
c_8 \in \mathbb{R}$ with $c_6 \notin \{-\tfrac{1}{2},0\}$ and
$c_2=c_5=0$ in
(\ref{Parameter-Ord22-A1-B1})--(\ref{Parameter-Ord22-A2-B2}) which
leads to
\begin{equation} \label{Parameter-Coeff-2.2.3}
    \alpha^T = \begin{bmatrix} 1 && c_6 && -c_6 \end{bmatrix},
    \quad
    A^{(0)} = \begin{bmatrix} 0 && 0 && 0 \\ c_7 && 0 && 0 \\
    \frac{1-2 c_6 c_7}{-2 c_6}-c_8 && c_8 && 0 \end{bmatrix},
    \quad
    B^{(0)} =
    \begin{bmatrix} 0 && 0 && 0 \\ \frac{c_1}{2}(1+\frac{1}{2 c_6}) && 0 && 0 \\
    \frac{c_1}{2}(1-\frac{1}{2 c_6}) && 0 && 0 \end{bmatrix} .
\end{equation}
\subsection{Coefficients for SRK Methods of Order (3,2)}
\label{Sec:Coeff-SRK-Ord-3-2}
If we consider the classification of the coefficients for explicit
SRK methods, we can see from Section~\ref{Sec:Coeff-SRK-Ord-2-2}
that in some of the resulting cases there are still degrees of
freedom in choosing the coefficients for $\alpha$ and $A^{(0)}$.
Therefore, we analyze now the classification for explicit SRK
methods (\ref{SRK-method-Ito-Wm-allg01}) with $s=3$ stages of order
$p_D=3$ and $p_S=2$. Thus, we additionally have to take into account
the well known deterministic order 3
conditions~\cite{Butcher87,HNW93}
\begin{alignat}{3}
    \alpha^T (A^{(0)} e)^2 &= \frac{1}{3} \, ,
    \label{cond-det-ord-3a} \\
    \alpha^T (A^{(0)} (A^{(0)} e)) &= \frac{1}{6} \, .
    \label{cond-det-ord-3b}
\end{alignat}
Clearly, these conditions can not be fulfilled in
case~(\ref{Parameter-Fall1}) where $\alpha_3=0$ as well as in
case~(\ref{Parameter-Fall2.1.1}) due to $A^{(0)}_{21}=0$ and in case
(\ref{Parameter-Fall2.2.2}) due to the restrictions for
$\alpha$ and $A^{(0)}$. However, in the case of parameter family
(\ref{Parameter-Fall2.1.2}) we obtain from (\ref{cond-det-ord-3a})
and (\ref{cond-det-ord-3b}) an SRK method of order (3,2) if in
(\ref{Parameter-Coeff-2.1.2}) it holds
\begin{equation}
    c_7 = \frac{1}{2} c_6 \pm \frac{1}{6} \sqrt{9 c_6^2-36 c_6
    +24} - \frac{1}{3 c_6} \, , \quad \quad c_8 = \frac{1}{3 c_6}
    \, .
\end{equation}
\quad \\
For the parameter families in case (\ref{Parameter-Fall2.2.1})
and (\ref{Parameter-Fall2.2.3}) we have to distinguish the
following three possibilities due to condition 8.\ of
Theorem~\ref{SRK-theorem-ito-ord2-Wm-main1} and due to
(\ref{cond-det-ord-3a}):
\begin{enumerate}[a)]
    \item \label{case-a)} $\alpha_2 = \frac{3}{4}$, $A^{(0)}_{21}=\frac{2}{3}$,
    $A^{(0)}_{31} = -A^{(0)}_{32}$. \\
    \item \label{case-b)} $\alpha_3=\frac{3}{4}-\alpha_2$,
    $A^{(0)}_{21}=\frac{2}{3}$,
    $A^{(0)}_{31}=\frac{2}{3}-A^{(0)}_{32}$. \\
    \item \label{case-c)} $\alpha_2 = \frac{1}{6} \frac{2-3(A^{(0)}_{31} +
    A^{(0)}_{32})}{A^{(0)}_{21} (A^{(0)}_{21}-A^{(0)}_{31}
    -A^{(0)}_{32})}$, $\alpha_3 = \frac{1}{6} \frac{3
    A^{(0)}_{21}-2}{(A^{(0)}_{31}+A^{(0)}_{32}) (A^{(0)}_{21} -
    A^{(0)}_{31}-A^{(0)}_{32})}$ if $A^{(0)}_{31}+A^{(0)}_{32} \neq
    A^{(0)}_{21} \neq 0$ and $A^{(0)}_{31}+A^{(0)}_{32} \neq 0$.
\end{enumerate}
If we consider now the case of the parameter family
(\ref{Parameter-Fall2.2.1}) then the conditions
(\ref{cond-det-ord-3a}) and (\ref{cond-det-ord-3b}) are fulfilled
for \ref{case-a)}) if
\begin{equation}
    c_6 = \frac{3}{4} \, , \quad \quad c_7 = \frac{1}{4 c_9} \, ,
    \quad \quad c_8 = \frac{2}{3} \, ,
\end{equation}
with $c_7 \notin \{-\tfrac{3}{4}, 0,\tfrac12
\} \cup \, ]-\tfrac{1}{4},0[ \,$ in (\ref{Parameter-Coeff-2.2.1}).
Further, the conditions (\ref{cond-det-ord-3a}) and
(\ref{cond-det-ord-3b}) are also fulfilled in the case
(\ref{Parameter-Fall2.2.1}) combined with \ref{case-b)}) if
\begin{equation}
    c_6 = \frac{3}{4} - \frac{1}{4 c_9} \, , \quad \quad c_7 =
    \frac{1}{4 c_9} \, , \quad \quad c_8 = \frac{2}{3} \, ,
\end{equation}
with $c_9 \neq 0$ and for $c_6 \in \, ]0,\tfrac{1}{4}[ \, \cup \,
]\tfrac{1}{4},\tfrac{3}{4}[\,$ in (\ref{Parameter-Coeff-2.2.1}).
Finally, the considered parameter family (\ref{Parameter-Fall2.2.1})
fulfills the conditions (\ref{cond-det-ord-3a}) and
(\ref{cond-det-ord-3b}) also for the case \ref{case-c)}) if
\begin{equation}\label{2.2.1.c}
    c_6 = \frac{1}{6} \frac{2-3 \lambda}{c_8 (c_8- \lambda)} \, , \quad
    \quad c_7 = \frac{1}{6} \frac{3 c_8 -2}{\lambda (c_8-\lambda)} \, ,\quad
    \quad c_9 = \frac{\lambda(c_8-\lambda)}{(3c_8-2)c_8}
\end{equation}
in (\ref{Parameter-Coeff-2.2.1}) for $\lambda \in \mathbb{R}$ with
$\lambda \notin \{0, \tfrac{2}{3}, c_8, \tfrac{2}{3}-c_8\}$,
$(\lambda-1)c_8\neq\lambda^2-\frac23$, $c_8 \notin \{0,
\tfrac{2}{3}\}$ and with $\lambda<\frac23$ if $c_8=1$;
$\tfrac{3c_8-2}{3(c_8-1)}\leq\lambda<\frac23$ if $\tfrac{2}{3} < c_8 < 1$ holds;
$\frac23<\lambda$ or $\lambda\leq\tfrac{3c_8-2}{3(c_8-1)}$ if $0 < c_8 <
\tfrac{2}{3}$ holds and with $\lambda<\frac23$ or $\lambda\geq\tfrac{3c_8-2}{3(c_8-1)}$
if $c_8 < 0$ or $c_8 > 1$ holds. Note that $\lambda = \frac{1-2 c_6
c_8}{2 c_7}$ is thus automatically
fulfilled in (\ref{Parameter-Coeff-2.2.1}). \\ \\
Finally, we consider the case of parameter family
(\ref{Parameter-Fall2.2.3}) which fulfills the additional order
three conditions (\ref{cond-det-ord-3a}) and (\ref{cond-det-ord-3b})
 in case \ref{case-a)}) for $c_6 = \frac{3}{4}$ and $c_7 =
\frac{2}{3}$ as well as in case \ref{case-c)}) with $c_6 =
\frac{1}{4 c_7 - \frac{4}{3}}$ if $c_7 \notin \{-\frac{1}{6},0,
\frac{1}{3}\}$ in (\ref{Parameter-Coeff-2.2.3}). For case
\ref{case-b)}) there exists no solution.
\section{Optimal SRK Schemes} \label{Sec:Optimal-schemes}
In the present section, coefficients for the SRK
method~(\ref{SRK-method-Ito-Wm-allg01}) of different orders of
convergence are presented. Due to some degrees of freedom in
choosing the coefficients, we consider additional conditions in
order to specify the free parameters of the SRK scheme. Clearly,
these additional conditions need not necessarily be fulfilled for
the desired order of convergence. However, coefficients fulfilling
also higher order conditions yield SRK schemes with the objective to
obtain smaller error constants and we call them optimal SRK schemes
in the following.
\subsection{Coefficients for Optimal SRK Schemes of Order (2,1)}
For SRK methods of order $p_D=2$ and $p_S=1$, we need 2 stages for
the drift part, however only one stage is needed for the
diffusion. Therefore, we consider only the case of $c_4=
\ldots=c_{11}=0$ in (\ref{Parameter-Ord21-all}).
Next, we want to specify $c_2$ and $c_3$. Therefore, we consider
additional order conditions which need not necessarily be fulfilled
for order $(2,1)$ schemes. Taking into account condition 9.\
yields $c_2 = c_3^2$. From condition 10.\ it follows that
$c_2 = c_1 \, c_3$. Further, one can consider the deterministic
order 3 conditions (\ref{cond-det-ord-3a}) and
(\ref{cond-det-ord-3b}) \cite{Butcher87,HNW93}, whereas only
(\ref{cond-det-ord-3a}) can be fulfilled which yields $c_2 =
\frac{2}{3}$. However, one can only combine two of the mentioned
additional conditions. Condition 9.\ together with 10.\ yields
$c_2=1$ and $c_3=c_1$, condition 9.\ together with
(\ref{cond-det-ord-3a}) yields $c_3 = \pm \sqrt{\frac{2}{3}}$ while
condition 10.\ together with (\ref{cond-det-ord-3a}) results
in $c_3=c_1 \frac{2}{3}$. One can easily verify that all
 the order 2 conditions 8.--50.\ are fulfilled with the
exception of conditions 11. and 13.--16.\ and condition 9.\
or 10.\ if (\ref{cond-det-ord-3a}) is fulfilled in combination
with only one of them.
Therefore, we consider the additional
condition~(\ref{cond-det-ord-3a}) which is fulfilled for $c_1 = 1$
and $c_2 = c_3 = \frac{2}{3}$. This leads to the SRK scheme
RDI1WM presented in Table~\ref{Table-Coeff-RDI1}, which is an
improved Euler-Maruyama scheme with two evaluations of the drift and
one of the diffusion coefficients for each step. Thus, it may be
superior to the widely used Euler-Maruyama scheme, especially in
practical applications where small noise is inherent to the system.
\begin{table}[tbp]
\renewcommand{\arraystretch}{1.3}
\begin{equation*}
\begin{array}{r|ccc|ccc|ccc}
    0 &&&&&&&&& \\
    \frac{2}{3} & \frac{2}{3} & & & \frac{2}{3} &&&&& \\
    \hline
    0 &&&&&&&&& \\
    0 & 0 & &  & 0 &&&&& \\
    \hline
    0 &&&&&&&&& \\
    0 & 0 &&& 0 &&&&& \\
    \hline
    & \frac{1}{4} & & \frac{3}{4} &
    1 & & 0 & 0 & \quad & 0 \\
    \cline{2-10}
    & & & & 0 & & 0 & 0 & & 0
\end{array}
\end{equation*}
\caption{Coefficients of the optimal SRK scheme RDI1WM with $p_D=2$
and $p_S=1$.} \label{Table-Coeff-RDI1}
\end{table}
\subsection{Coefficients for Optimal SRK Schemes of Order (2,2) and (3,2)}
\begin{figure}[tbp]
\begin{center}
\includegraphics[width=6.8cm]{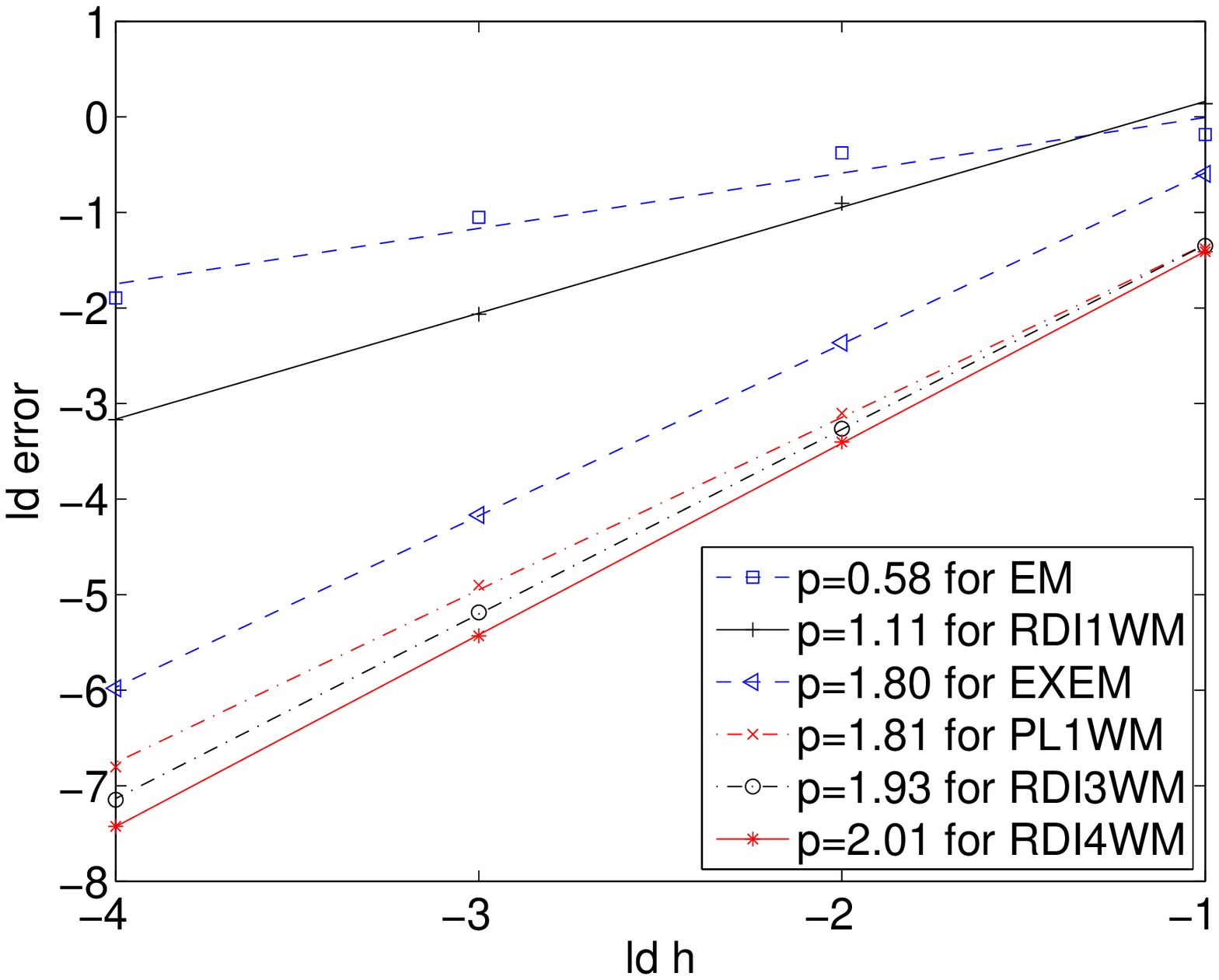}
\includegraphics[width=6.8cm]{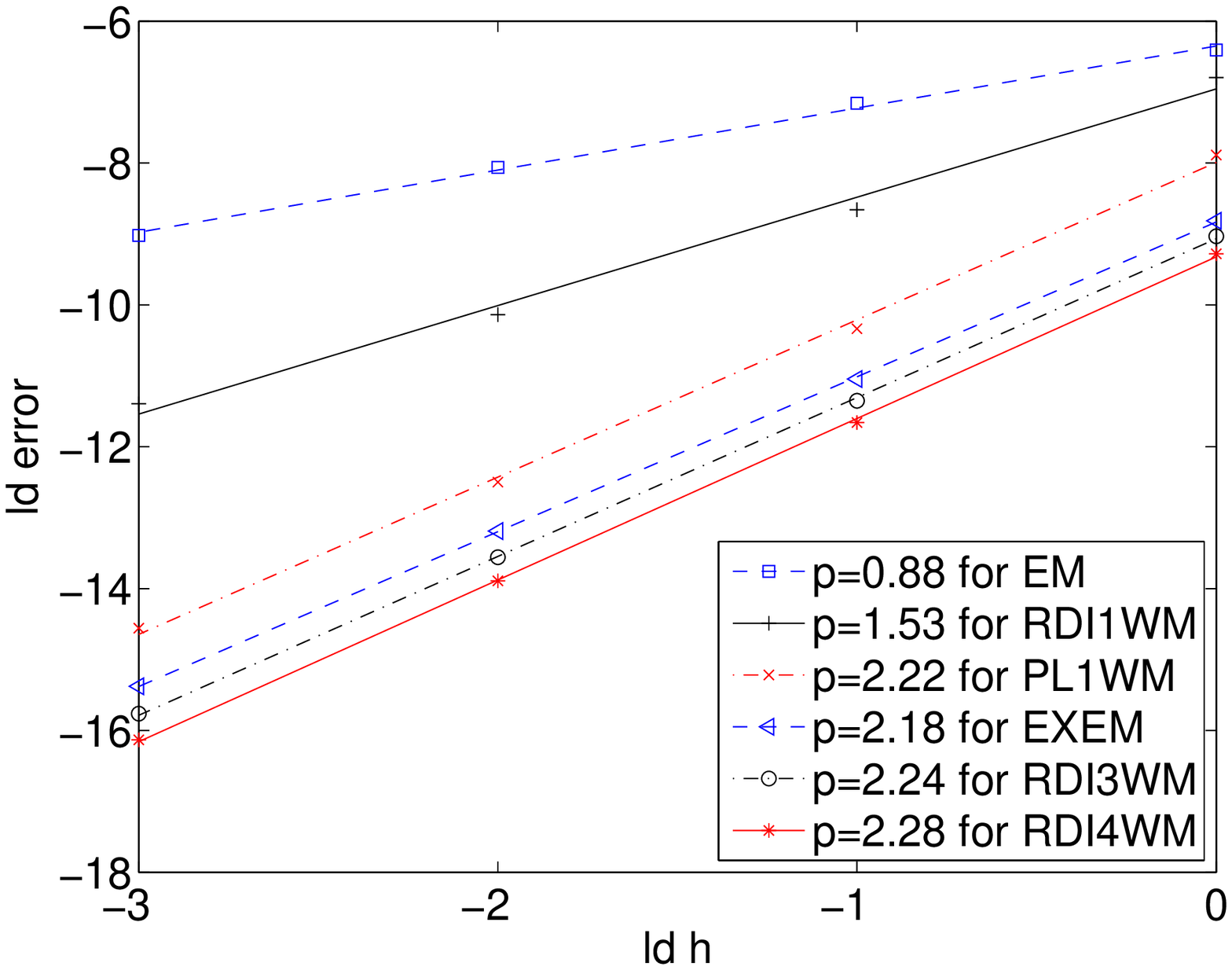}
\caption{Orders of convergence for
SDE~(\ref{Simu:nonlinear-SDE2}) and SDE~(\ref{Simu:dm-SDE}).}
\label{Bild001}
\end{center}
\end{figure}

If we consider the order three tree $(\sigma_{j_1}, \sigma_{j_2},
\{\tau, \sigma_{j_4} \}_{j_3})$ (see \cite{Roe06a,Roe03} for
details) in the case of $j_1=j_2=j_3=j_4$, then we obtain the
corresponding order condition
\begin{equation} \label{tree-cond-1}
    {\beta^{(2)}}^T
    ((A^{(1)} e) (B^{(1)} e))
    ({\beta^{(1)}}^T e)^2 = \frac{2}{3} .
\end{equation}
For the coefficient
families~(\ref{Parameter-Coeff-1})-(\ref{Parameter-Coeff-2.2.3}),
this order condition is fulfilled if $c_3 = \pm \sqrt{\frac{2}{3}}$.
For the tree $(\sigma_{j_1},\{\sigma_{j_2},\sigma_{j_2},
\sigma_{j_3}, \sigma_{j_3} \}_{j_1})$ (see \cite{Roe06a,Roe03}) we
calculate in the case of $j_1 \neq j_2$ and $j_2=j_3$ the
following order three condition
\begin{equation} \label{tree-cond-2}
    ({\beta^{(1)}}^T e) ({\beta^{(3)}}^T
    ( B^{(2)} e )^4) = 1
\end{equation}
which is fulfilled if $c_4 = \pm \sqrt{2}$.
Due to some symmetry in the SRK schemes, we obtain always the same
SRK schemes regardless what sign we choose for $c_3$, $c_4$ and
$c_1$. In the following, we choose $c_3 = \sqrt{\frac{2}{3}}$, $c_4
= \sqrt{2}$ and $c_1=1$.\\\\
Then, the parameter family (\ref{Parameter-Coeff-1}) definitely
provides the optimal SRK scheme {RDI2WM of order $p_D=p_S=2$
presented in Table~\ref{Table-Coeff-RDI2}.
\begin{table}[tbp]
\renewcommand{\arraystretch}{1.3}
\begin{equation*}
\begin{array}{r|ccccc|ccccc|cccccc}
    0 & & & &&  &&&&&  &&& \\
    1 & 1 & &&&  & 1 & &&  &&& \\
    0 & 0 & & 0 &&  & 0 & & 0 &  &&& \\
    \hline
    0 &&&  &&&&&  &&& \\
    \frac{2}{3} & \frac{2}{3} & & & &  & \sqrt{\frac{2}{3}} &&&&  &&& \\
    \frac{2}{3} & \frac{2}{3} & & 0 & &  & -\sqrt{\frac{2}{3}} && 0 &&  &&& \\
    \hline
    0 & & &&&  &&&&&  &&& \\
    0 & 0 & &&&  & \sqrt{2} &&&&  &&& \\
    0 & 0 & & 0 & & & -\sqrt{2} & & 0 & & &&& \\
    \hline
    & \frac{1}{2} & & \frac{1}{2} & & 0 &
    \frac{1}{4} & & \frac{3}{8} & & \frac{3}{8} & & 0 & &
    \frac{\sqrt{6}}{4} & & -\frac{\sqrt{6}}{4} \\
    \cline{2-17}
    & & & & & & -\frac{1}{4} & & \frac{1}{8} & & \frac{1}{8}
    & & 0 & & \frac{\sqrt{2}}{4} & & -\frac{\sqrt{2}}{4}
\end{array}
\end{equation*}
\caption{Coefficients of the optimal SRK scheme RDI2WM with $p_D=2$
and $p_S=2$.} \label{Table-Coeff-RDI2}
\end{table}
\\ \\
Next, we calculate SRK schemes of order $p_D=3$ and $p_S=2$ for
the family~(\ref{Parameter-Coeff-2.2.1}) in the case
(\ref{2.2.1.c}). Again, we try to specify the remaining
coefficients in the deterministic part of the scheme by
additionally considering the order four
conditions~\cite{Butcher87,HNW93}
\begin{alignat}{3}
    \alpha^T (A^{(0)} (A^{(0)}e)^2) &=
    \frac{1}{12} \, , \label{Zusatz-Gl-001a} \\
    \alpha^T ((A^{(0)}e)(A^{(0)}
    (A^{(0)}e) ) ) &= \frac{1}{8} \, . \label{Zusatz-Gl-001}
\end{alignat}
These conditions are fulfilled if $\lambda=\tfrac{3}{4}$ and
$c_8=\tfrac{1}{2}$. As a result of this, we obtain the coefficients
of the SRK scheme RDI3WM presented in Table~\ref{Table-Coeff-RDI3}.
\begin{table}[tbp]
\renewcommand{\arraystretch}{1.3}
\begin{equation*}
\begin{array}{r|ccccc|ccccc|cccccc}
    0 & & & &&  &&&&&  &&& \\
    \frac{1}{2} & \frac{1}{2} & &&&  & \frac{9-2\sqrt{15}}{14} & &&  &&& \\
    \frac{3}{4} & 0 & & \frac{3}{4} &&  & \frac{18+3 \sqrt{15}}{28} & & 0 &  &&& \\
    \hline
    0 &&&  &&&&&  &&& \\
    \frac{2}{3} & \frac{2}{3} & & & &  & \sqrt{\frac{2}{3}} &&&&  &&& \\
    \frac{2}{3} & \frac{2}{3} & & 0 & &  & -\sqrt{\frac{2}{3}} && 0 &&  &&& \\
    \hline
    0 & & &&&  &&&&&  &&& \\
    0 & 0 & &&&  & \sqrt{2} &&&&  &&& \\
    0 & 0 & & 0 & & & -\sqrt{2} & & 0 & & &&& \\
    \hline
    & \frac{2}{9} & & \frac{1}{3} & & \frac{4}{9} &
    \frac{1}{4} & & \frac{3}{8} & & \frac{3}{8}
    & & 0 & & \frac{\sqrt{6}}{4} & & -\frac{\sqrt{6}}{4} \\
    \cline{2-17}
    & & & & & & -\frac{1}{4} & & \frac{1}{8} & & \frac{1}{8}
    & & 0 & & \frac{\sqrt{2}}{4} & & -\frac{\sqrt{2}}{4}
\end{array}
\end{equation*}
\caption{Coefficients of the optimal SRK scheme RDI3WM with $p_D=3$
and $p_S=2$.} \label{Table-Coeff-RDI3}
\end{table}
\\ \\
However, if we claim for the family~(\ref{Parameter-Coeff-2.2.1}) in
the case of (\ref{2.2.1.c}) that the order four
conditions~(\ref{Zusatz-Gl-001a}) and
\begin{equation} \label{Zusatz-Gl-002}
    \begin{split}
    \alpha^T (A^{(0)}e)^3 &= \frac{1}{4}
    \end{split}
\end{equation}
are fulfilled instead of~(\ref{Zusatz-Gl-001}), then we get the
coefficients $\lambda=1$ and $c_8=\frac{1}{2}$.
As a result of this, we obtain the coefficients of the SRK scheme
RDI4WM presented in Table~\ref{Table-Coeff-RDI4}. Here, the
deterministic part of scheme RDI4WM coincides with the well
known Simpson scheme for ODEs~\cite{HNW93}.
\begin{table}[tbp]
\renewcommand{\arraystretch}{1.3}
\begin{equation*}
\begin{array}{r|ccccc|ccccc|cccccc}
    0 & & & &&  &&&&&  &&& \\
    \frac{1}{2} & \frac{1}{2} & &&&  & \frac{6-\sqrt{6}}{10} & &&  &&& \\
    1 & -1 & 2 &  &&  & \frac{3+2\sqrt{6}}{5} & & 0 &  &&& \\
    \hline
    0 &&&  &&&&&  &&& \\
    \frac{2}{3} & \frac{2}{3} & & & &  & \sqrt{\frac{2}{3}} &&&&  &&& \\
    \frac{2}{3} & \frac{2}{3} & & 0 & &  & -\sqrt{\frac{2}{3}} && 0 &&  &&& \\
    \hline
    0 & & &&&  &&&&&  &&& \\
    0 & 0 & &&&  & \sqrt{2} &&&&  &&& \\
    0 & 0 & & 0 & & & -\sqrt{2} & & 0 & & &&& \\
    \hline
    & \frac{1}{6} & & \frac{2}{3} & & \frac{1}{6} &
    \frac{1}{4} & & \frac{3}{8} & & \frac{3}{8}
    & & 0 & & \frac{\sqrt{6}}{4} & & -\frac{\sqrt{6}}{4} \\
    \cline{2-17}
    & & & & & & -\frac{1}{4} & & \frac{1}{8}
    & & \frac{1}{8} & & 0 & & \frac{\sqrt{2}}{4} & & -\frac{\sqrt{2}}{4}
\end{array}
\end{equation*}
\caption{Coefficients of the optimal SRK scheme RDI4WM with $p_D=3$
and $p_S=2$.} \label{Table-Coeff-RDI4}
\end{table}
\section{Numerical example} \label{Sec:Numerical-Example}
In the following, some of the SRK schemes presented in
Section~\ref{Sec:Optimal-schemes} are applied to test equations in
order to analyze their order of convergence in comparison to
some well known schemes.
Therefore, the functional $u = \E(f(X_t))$ is approximated by a
Monte Carlo simulation based on the optimal SRK schemes RDI1WM of
order 1 and RDI3WM and RDI4WM of order 2. The optimal SRK schemes
are compared to the second order SRK scheme PL1WM due to
Platen~\cite{KP99}, the Euler--Maruyama scheme EM of order 1 and the
ex\-tra\-po\-la\-ted Euler-Maruyama scheme EXEM \cite{KP99} also
attaining order 2. The SRK scheme PL1WM is contained in the class of
SRK schemes (\ref{SRK-method-Ito-Wm-allg01}) with coefficients $c_1
= c_3 = c_4 = 1$ in (\ref{Parameter-Coeff-1}) due to case
(\ref{Parameter-Fall1}). The extra\-po\-lated Euler-Maruyama
approximation is given by $2 \E(f(Z_t^{h/2}))- \E(f(Z_t^{h}))$ based
on the Euler-Maruyama approximations $Z_t^{h/2}$ and $Z_t^h$
calculated with step sizes $h$ and $h/2$.} The sample average
$u_{M,h} = \frac{1}{M} \sum_{k=1}^M f(Y_t(\omega_k))$, $\omega_k \in
\Omega$, of $M$ independent simulated realizations of the considered
approximation $Y_t$ is calculated in order to estimate the
expectation.
In the following, we denote by $\hat{\mu} = u_{M,h} - \E(f(X_t))$
the mean error and by $\hat{\sigma}^2_{\mu}$ the empirical
variance of the mean error. Further, we calculate the confidence
interval with boundaries $a$ and $b$ to the level of 90\% for the
estimated error $\hat{\mu}$ (see \cite{KP99,Roe03} for
details). \\ \\
As the first example, we consider the non-linear
SDE~\cite{KP99,MacNav01,Roe06b}
\begin{equation} \label{Simu:nonlinear-SDE2}
    dX_t = \left( \tfrac{1}{2} X_t + \sqrt{X_t^2 + 1} \right) \,
    dt + \sqrt{X_t^2 + 1} \, dW_t, \qquad X_0=0,
\end{equation}
on the time interval $I=[0,2]$ with the solution $X_t = \sinh (t +
W_t)$. Here, we choose $f(x)=p(\arsinh(x))$, where $p(z) = z^3 -
6z^2 + 8z$ is a polynomial. Then the expectation of the solution
can be calculated as
\begin{equation}
    \E(f(X_t)) = t^3 - 3t^2 + 2t \,\,.
\end{equation}
The solution $\E(f(X_t))$ is approximated with step sizes $2^{-1},
\ldots, 2^{-4}$ and $M=10^9$ simulations are performed in order to
determine the systematic error of the considered schemes at time
$t=2$. The results for the applied schemes are presented in
Table~\ref{Table1}. The orders of convergence correspond to the
slope of the regression lines plotted in Figure~\ref{Bild001}
where we get the order $0.58$ for the EM scheme, order $1.11$ for
RDI1WM, order $1.80$ for EXEM, order $1.81$ for PL1WM, order $1.93$
for RDI3WM and order $2.01$ for the scheme RDI4WM.
\begin{table}
\caption{Mean errors, empirical variances and confidence intervals
for SDE~(\ref{Simu:nonlinear-SDE2}).} \label{Table1}
\begin{center}
\setlength{\extrarowheight}{-3pt}
\begin{tabular}{c|c|c|c|c|c}
   & $h$ & $|\hat{\mu}|$ & $\hat{\sigma}_{\mu}^2$  & $a$ & $b$ \\
   \hline
\multirow{4}{*}{EM}
   & $2^{-1}$ & 8.797E-01 & 6.534E-07 & -8.799E-01 & -8.795E-01 \\
   & $2^{-2}$ & 7.705E-01 & 1.592E-06 & -7.708E-01 & -7.702E-01 \\
   & $2^{-3}$ & 4.825E-01 & 1.599E-06 & -4.828E-01 & -4.822E-01 \\
   & $2^{-4}$ & 2.691E-01 & 1.754E-06 & -2.694E-01 & -2.688E-01 \\
    \hline
\multirow{4}{*}{RDI1WM}
   & $2^{-1}$ & 1.101E-00 & 1.381E-06 & -1.101E-00 & -1.100E-00\\
   & $2^{-2}$ & 5.342E-01 & 2.080E-06 & -5.346E-01 & -5.339E-01\\
   & $2^{-3}$ & 2.390E-01 & 3.297E-06 & -2.394E-01 & -2.386E-01\\
   & $2^{-4}$ & 1.112E-01 & 2.984E-06 & -1.116E-01 & -1.107E-01\\
    \hline
\multirow{4}{*}{EXEM}
   & $2^{-1}$ & 1.359E-00 & 2.990E-06 & -1.359E-00 & -1.359E-00\\
   & $2^{-2}$ & 6.614E-01 & 7.315E-06 & -6.620E-01 & -6.607E-01\\
   & $2^{-3}$ & 1.945E-01 & 8.629E-06 & -1.952E-01 & -1.938E-01\\
   & $2^{-4}$ & 5.570E-02 & 9.014E-06 & -5.641E-02 & -5.499E-02\\
    \hline
\multirow{4}{*}{PL1WM}
   & $2^{-1}$ & 3.837E-01 & 1.885E-06 & -3.841E-01 & -3.834E-01\\
   & $2^{-2}$ & 1.165E-01 & 3.207E-06 & -1.169E-01 & -1.161E-01\\
   & $2^{-3}$ & 3.348E-02 & 2.475E-06 & -3.386E-02 & -3.311E-02\\
   & $2^{-4}$ & 8.949E-03 & 3.447E-06 & -9.390E-03 & -8.509E-03\\
    \hline
\multirow{4}{*}{RDI3WM}
   & $2^{-1}$ & 3.926E-01 & 1.400E-06 & -3.929E-01 & -3.923E-01 \\
   & $2^{-2}$ & 1.041E-01 & 2.787E-06 & -1.045E-01 & -1.037E-01 \\
   & $2^{-3}$ & 2.748E-02 & 2.427E-06 & -2.785E-02 & -2.711E-02 \\
   & $2^{-4}$ & 7.054E-03 & 1.813E-06 & -7.373E-03 & -6.734E-03\\
    \hline
\multirow{4}{*}{RDI4WM}
   & $2^{-1}$ & 3.760E-01 & 1.488E-06 & -3.762E-01 & -3.757E-01 \\
   & $2^{-2}$ & 9.454E-02 & 2.823E-06 & -9.494E-02 & -9.414E-02 \\
   & $2^{-3}$ & 2.318E-02 & 2.441E-06 & -2.355E-02 & -2.281E-02 \\
   & $2^{-4}$ & 5.816E-03 & 1.816E-06 & -6.135E-03 & -5.496E-03
\end{tabular}
\end{center}
\end{table}
\\ \\
As a second example, a multi-dimensional SDE with a
$2$-dimensional dri\-ving Wiener process is considered:
\begin{equation} \label{Simu:dm-SDE}
    \begin{split}
    d \begin{pmatrix} X_t^1 \\ X_t^2 \end{pmatrix}
    &= \begin{pmatrix} -\frac{273}{512} & 0 \\
    -\frac{1}{160} & -\frac{785}{512}+\frac{\sqrt{2}}{8} \end{pmatrix} \,
    \begin{pmatrix} X_t^1 \\ X_t^2 \end{pmatrix} \, dt +
    \begin{pmatrix} \frac{1}{4} X_t^1 & \frac{1}{16} X_t^1 \\
    \frac{1-2\sqrt{2}}{4} X_t^2 & \frac{1}{10} X_t^1 + \frac{1}{16} X_t^2
    \end{pmatrix} \, d \begin{pmatrix} W_t^1 \\ W_t^2
    \end{pmatrix}
    \end{split}
\end{equation}
with initial value $X_0=(1,1)^T$. This SDE system is of special
interest due to the fact that it has no commutative noise. Here,
we are interested in the second moments which depend on both, the
drift and the diffusion function (see \cite{KP99} for details).
Therefore, we choose $f(x) = (x^1)^2$
and obtain
\begin{equation}
    \begin{split}
    \E(f(X_t)) = \exp(- t) \, .
    \end{split}
\end{equation}
We approximate $\E(f(X_t))$ at $t=4$ by $M=5\cdot 10^8$ simulated
trajectories with step sizes $2^{-0}, \ldots, 2^{-3}$. The results
for the considered schemes are presented in Table~\ref{Table2}
and Figure~\ref{Bild001}. Here, the order of convergence is $0.88$
for the Euler-Maruyama scheme, $1.53$ for RDI1WM, $2.22$ for PL1WM,
$2.18$ for EXEM, $2.24$ for RDI3WM and order $2.28$ for the
optimal SRK scheme RDI4WM.
\begin{table}
\caption{Mean errors, empirical variances and confidence intervals
for SDE~(\ref{Simu:dm-SDE}).} \label{Table2}
\begin{center}
\setlength{\extrarowheight}{-3pt}
\begin{tabular}{c|c|c|c|c|c}
   & $h$ & $|\hat{\mu}|$ & $\hat{\sigma}_{\mu}^2$  & $a$ & $b$ \\
   \hline
\multirow{4}{*}{EM}
   & $2^{-0}$ & 1.178E-02 & 3.946E-11 & -1.178E-02 & -1.178E-02 \\
   & $2^{-1}$ & 7.002E-03 & 6.669E-11 & -7.004E-03 & -7.000E-03 \\
   & $2^{-2}$ & 3.738E-03 & 5.799E-11 & -3.740E-03 & -3.736E-03 \\
   & $2^{-3}$ & 1.922E-03 & 8.614E-11 & -1.925E-03 & -1.920E-03 \\
    \hline
\multirow{4}{*}{RDI1WM}
   & $2^{-0}$ & 9.000E-03 & 1.275E-10 & 8.998E-03 & 9.004E-03\\
   & $2^{-1}$ & 2.472E-03 & 1.127E-10 & 2.470E-03 & 2.475E-03 \\
   & $2^{-2}$ & 8.870E-04 & 8.278E-11 & 8.848E-04 & 8.891E-04\\
   & $2^{-3}$ & 3.714E-04 & 8.926E-11 & 3.691E-04 & 3.736E-04\\
    \hline
\multirow{4}{*}{EXEM}
   & $2^{-0}$ & 2.223E-03 & 2.871E-10 & -2.227E-03 & -2.219E-03\\
   & $2^{-1}$ & 4.733E-04 & 2.500E-10 & -4.771E-04 & -4.696E-04\\
   & $2^{-2}$ & 1.071E-04 & 3.585E-10 & -1.116E-04 & -1.026E-04\\
   & $2^{-3}$ & 2.348E-05 & 3.919E-10 & -2.818E-05 & -1.879E-05\\
    \hline
\multirow{4}{*}{PL1W1}
   & $2^{-0}$ & 4.230E-03 & 6.967E-11 & 4.228E-03 & 4.232E-03 \\
   & $2^{-1}$ & 7.736E-04 & 8.594E-11 & 7.714E-04 & 7.758E-04 \\
   & $2^{-2}$ & 1.728E-04 & 8.412E-11 & 1.706E-04 & 1.750E-04 \\
   & $2^{-3}$ & 4.148E-05 & 8.356E-11 & 3.932E-05 & 4.365E-05 \\
    \hline
\multirow{4}{*}{RDI3WM}
   & $2^{-0}$ & 1.909E-03 & 3.700E-11 & -1.910E-03 & -1.907E-03\\
   & $2^{-1}$ & 3.822E-04 & 6.597E-11 & -3.841E-04 & -3.803E-04 \\
   & $2^{-2}$ & 8.282E-05 & 6.356E-11 & -8.471E-05 & -8.093E-05 \\
   & $2^{-3}$ & 1.797E-05 & 8.787E-11 & -2.019E-05 & -1.574E-05\\
    \hline
\multirow{4}{*}{RDI4WM}
   & $2^{-0}$ & 1.608E-03 & 4.285E-11 & -1.609E-03 & -1.606E-03 \\
   & $2^{-1}$ & 3.089E-04 & 6.812E-11 & -3.108E-04 & -3.069E-04 \\
   & $2^{-2}$ & 6.583E-05 & 6.403E-11 & -6.773E-05 & -6.394E-05 \\
   & $2^{-3}$ & 1.392E-05 & 8.803E-11 & -1.615E-05 & -1.170E-05 \\
\end{tabular}
\end{center}
\end{table}
\begin{figure}[tbp]
\begin{center}
\includegraphics[width=6.8cm]{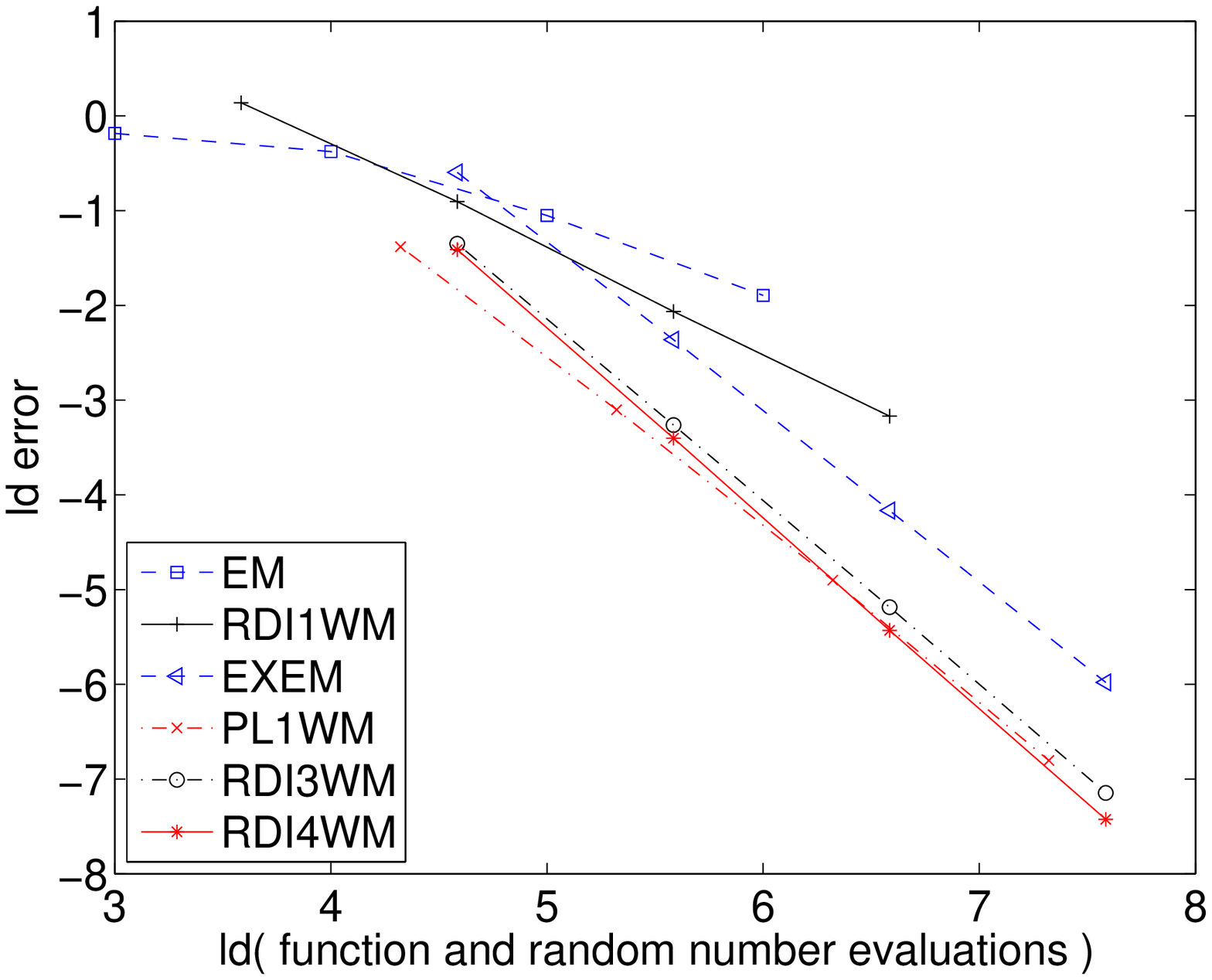}
\includegraphics[width=6.8cm]{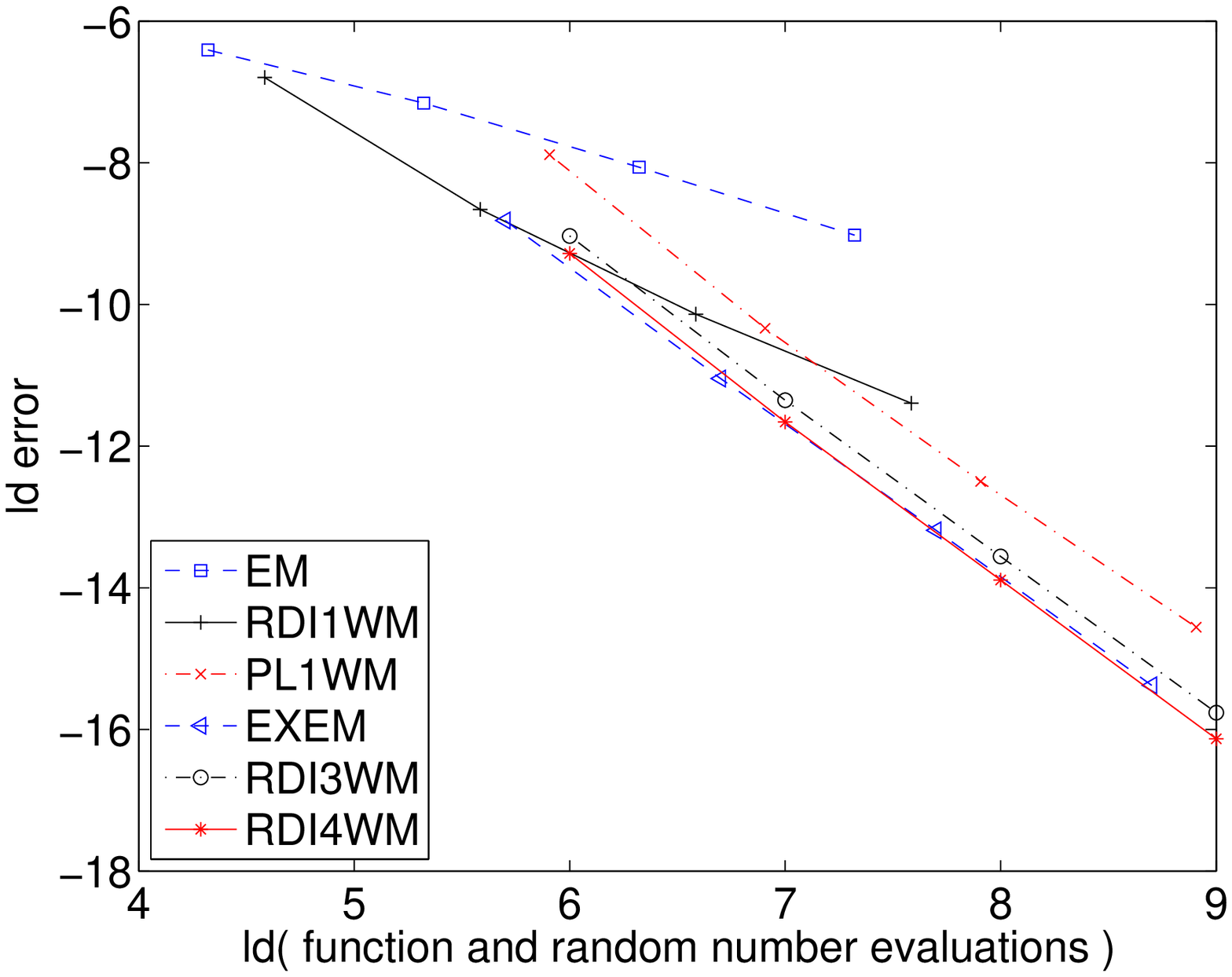}
\caption{Computational effort per simulation path versus precision for
SDE~(\ref{Simu:nonlinear-SDE2}) and SDE~(\ref{Simu:dm-SDE}).}
\label{Bild002}
\end{center}
\end{figure}
\quad \\ \\
Due to the results in Figure~\ref{Bild001}, we can see that for both
test equations the so called optimal SRK scheme RDI1WM attains much
better orders of convergence than the well known order one EM
scheme. The same holds for the optimal SRK schemes RDI3WM and
RDI4WM compared to the order two schemes EXEM and PL1WM. Clearly,
the optimal SRK schemes RDI1WM, RDI3WM and RDI4WM need some
additional computational effort compared to the schemes EM, EXEM and
PL1WM, respectively. Therefore, we take the number of evaluations of
the drift function $a$ and of each diffusion function $b^j$, $1 \leq
j \leq m$, as well as the number of random variables that have to be
simulated as a measure for the computational effort. Then we can
compare the computational effort versus the errors of the analyzed
schemes. The results are presented in Figure~\ref{Bild002}, and
again, RDI1WM performs much better than the scheme EM for both test
equations. Further, RDI3WM yields similar results like RDI4WM which
is for higher precisions slightly better than PL1WM and
significantly better than the scheme EXEM for the test
equation~(\ref{Simu:nonlinear-SDE2}). Considering the
multi-dimensional test equation (\ref{Simu:dm-SDE}), the scheme
RDI3WM is again close to RDI4WM which performs for higher precisions
slightly better than EXEM. However both optimal SRK schemes RDI3WM
and RDI4WM are significantly better than the SRK scheme PL1WM.
\section{Conclusion}
\label{Sec:Conclusion}
In the present work, a full classification of the coefficients for a
class of explicit SRK methods of order $(1,1)$ for $s=1$ and order
$(2,1)$ for $s=2$ stages as well as for the orders $(2,2)$ and
$(3,2)$ with $s=3$ stages is calculated. Based on this
classification, coefficients for so called optimal SRK schemes are
determined by considering additional higher order conditions.
Optimal coefficients for SRK methods of order $(2,1)$, $(2,2)$ and
$(3,2)$ are calculated and similarly to the deterministic
setting~\cite{HNW93}, better convergence results are expected for
these schemes in general. Finally, the SRK schemes based on optimal
coefficients are applied to some test equations. Here, it turned out
that the proposed optimal SRK schemes attain higher orders of
convergence than the well known schemes under consideration
and they also perform very well if the computational effort is taken
into account. \\ \\
For future research, it would be interesting to extend the presented
classification to diagonal or fully implicit SRK methods. Further,
the given classification may be applied in order to determine
coefficients for SRK methods with optimal stability properties.
\section*{Acknowledgements}
The authors are very grateful to the unknown referees for their
fruitful comments and suggestions.

\end{document}